\documentclass{ifacconf}

\usepackage{graphicx}      
\usepackage{natbib}        

\usepackage{amsmath}

\begin{document}
\begin{frontmatter}

\title{Flatness Analysis for the Sampled-data Model of a Single Mast Stacker Crane\thanksref{footnoteinfo}} 


\thanks[footnoteinfo]{This work has been supported by the Austrian Science Fund (FWF) under grant number P 32151.}

\author[First]{Johannes Diwold} 
\author[Second]{Bernd Kolar} 
\author[First]{Markus Schöberl}

\address[First]{Institute of Automatic Control and Control Systems Technology,	Johannes Kepler University Linz, Altenbergerstraße 66, 4040 Linz, Austria (e-mail: johannes.diwold@jku.at, markus.schoeberl@jku.at)}
\address[Second]{Magna Powertrain Engineering Center Steyr GmbH \& Co KG, Steyrer	Str. 32, 4300 St. Valentin, Austria (e-mail: bernd\_kolar@ifac-mail.org)}

\begin{abstract}                
We show that the Euler-discretization of the nonlinear continuous-time model of a single mast stacker crane is flat. The construction of the flat output is based on a transformation of a subsystem into the linear time-variant discrete-time controller canonical form. Based on the derived flat output, which is also a function of backward-shifts of the system variables, we discuss the planning of trajectories to achieve a transition between two rest positions and compute the corresponding discrete-time feedforward control.
\end{abstract}


\end{frontmatter}

\section{Introduction}
In the 1990s, the concept of flatness has been introduced by Fliess,
L\'{e}vine, Martin and Rouchon for nonlinear continuous-time systems
(see e.g. \cite{FliessLevineMartinRouchon:1995} and \cite{FliessLevineMartinRouchon:1999}).
A nonlinear continuous-time system
\begin{equation}
\dot{x}=f(x,u)\label
{eq:system_equation_continuous}
\end{equation}
with $\dim(x)=n$ states, $\dim(u)=m$ inputs and smooth functions
$f$ is flat, if there exists a one-to-one correspondence between
solutions $(x(t),u(t))$ of \eqref{eq:system_equation_continuous}
and solutions $y(t)$ of a trivial system (sufficiently smooth but
otherwise arbitrary trajectories). This one-to-one correspondence
implies that for flat systems there exist flat outputs 
\begin{equation}
\label{eq:cont_flat_output}
y=\varphi(x,u,\dot{u},\dots,u^{(q)})
\end{equation}
such that all states $x$ and inputs $u$ can be parameterized by
$y$ and its time derivatives, i.e., 
\begin{equation}
(x,u)=F(y,\dot{y},\dots,y^{(r)})\,.\label{eq:parameterizing_map_continuous}
\end{equation}
We call $F$ the parameterization w.r.t. the flat output
$y$. %

In the discrete-time framework one could define flatness simply by replacing the time-derivatives in \eqref{eq:cont_flat_output} and \eqref{eq:parameterizing_map_continuous} by forward-shifts, see e.g.  \cite{Sira-RamirezAgrawal:2004} and \cite{KaldmaeKotta:2013}. However, it has recently been pointed out in \cite{GuillotMillerioux2020} that this definition covers only a subclass of flat systems, which we also refer to as forward-flat systems in the sequel. Accordingly, in \cite{DiwoldKolarSchoeberl2022} we have derived a concise geometric definition of flatness, which is (like in the continuous-time case) solely based on the existence of a one-to-one correspondence between solutions of a flat system and solutions of a trivial system. The proposed approach leads to a definition of flatness which includes also backward-shifts of system variables in the flat output.

Both continuous-time and discrete-time flatness allow an elegant solution to motion planning problems as well as the systematic design of tracking controllers. For the practical implementation of a continuous-time controller, however, it is necessary to evaluate the control law at a sufficiently high sampling rate. If, nevertheless, the control loop is limited to lower sampling rates (possibly also due to an online-optimization), a discrete-time control law based on a suitable discretization might be preferable, see e.g. \cite{DiwoldKolarSchoeberl2022-1}. W.r.t. the discrete-time flatness-based controller design, the main challenge is to find a flat discretization as well as the corresponding flat output. So far, computationally efficient methods for the systematic construction of flat outputs have only been derived for forward-flat systems, see \cite{KolarSchoberlDiwold:2019} and \cite{KolarDiwoldSchoberl:2019}. However, for the considered sampled-data model of a single mast stacker crane (obtained by an explicit Euler-discretization) a forward-flat output does not exist. Nevertheless, the considered discrete-time system is flat in the extended sense including backward-shifts, as we show in this contribution. 

We consider the finite-dimensional, continuous-time model of a single mast stacker crane discussed in \cite{StaudeckerSchlacherHansl2008} and \cite{RamsSchoeberlSchlacher2018}. While in general the construction of flat outputs is a difficult task, two flat outputs are known for the continuous-time system. The first is a configuration-flat output, which can be derived using the methods of \cite{RathinamMurray1996} or \cite{GstoettnerKolarSchoeberl2021-2}. However, since the corresponding parameterizing map \eqref{eq:parameterizing_map_continuous} has singularities at stationary points, the configuration-flat output is not suitable for practical purposes. Rather suitable is the second known flat output, whose first component corresponds to the vertical position of the lifting unit. As shown in \cite{StaudeckerSchlacherHansl2008}, the second component is of the form $y^2=\varphi^2(x,u,\dot{u})$ and can be constructed systematically by transforming a subsystem into the time-variant controller canonical form. 
In this contribution we will transfer the ideas of \cite{StaudeckerSchlacherHansl2008} to the discrete-time framework, and construct a flat output for the sampled-data model of a stacker crane by means of the time-variant discrete-time controller canonical form. While the first component of the obtained flat output will be preserved, the second component will also be a function of backward-shifts of the system variables. To highlight the practical applicability of the derived flat output, we compute a feedforward control in order to transfer the single mast stacker crane between two rest positions.


The paper is organized as follows: In Section \ref{sec:model_stacker_crane} we recall the continuous-time model of a single mast stacker crane used in \cite{StaudeckerSchlacherHansl2008} and \cite{RamsSchoeberlSchlacher2018}, and recapitulate known results regarding the flatness of the system. Next, in Section \ref{sec:flatness_time_variant} we briefly recall and extend the notion of discrete-time flatness as proposed in \cite{DiwoldKolarSchoeberl2022} to time-variant nonlinear systems, and discuss the controller canonical form for linear time-variant discrete-time single-input systems. Based on the results of Section \ref{sec:flatness_time_variant}, we construct a flat output for the sampled-data model of a stacker crane in Section \ref{sec:sampled_data_model_stacker_crane_}. Finally, we compute a feedforward control, which transfers the single mast stacker crane between two rest positions.

\section{Model of a Single Mast Stacker Crane}
\label{sec:model_stacker_crane}
In this section, we recall the finite-dimensional, continuous-time model of a single mast stacker crane used in \cite{StaudeckerSchlacherHansl2008} and \cite{RamsSchoeberlSchlacher2018}. Furthermore, we recapitulate known results regarding the flatness of the system and show how a flat output can be derived by means of the controller canonical form for linear time-variant systems. In the remainder of the paper, we will transfer this idea to the discrete-time framework and show that the explicit Euler-discretization of the considered continuous-time model is flat as well.
\begin{figure}
	\begin{center}
		\includegraphics[width=8.4cm]{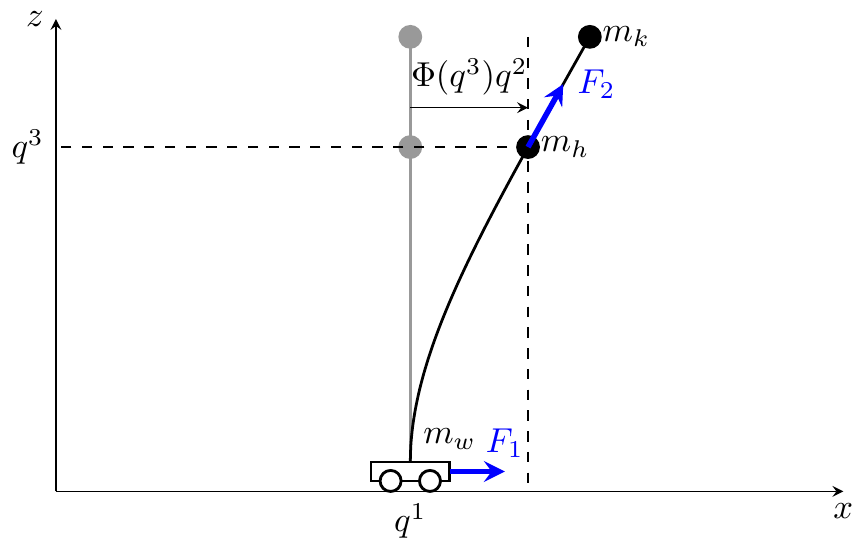}    
		\caption{Schematic of a single mast stacker crane.} 
		\label{fig:schematic}
	\end{center}
\end{figure}
\subsection{Derivation of the Continuous-time Model}
In Fig. \ref{fig:schematic} the setup of a single mast stacker crane is illustrated. With $q^1$ we denote the position of the rigid driving unit with mass $m_w$, while $q^3$ describes the vertical position of the lifting unit with mass $m_h$. The mast with length $L$ is modeled as an Euler-Bernoulli beam with the mass density $\rho A$ and the flexural rigidity $EI$. Like in \cite{StaudeckerSchlacherHansl2008} we approximate the deformation of the beam by a first-order Rayleigh-Ritz ansatz
\begin{equation*}
w(z,t)=\Phi(z)q^2(t)\,.
\end{equation*}
Evaluating the Euler-Lagrange equations for the finite-dimensional approximation yields a finite-dimensional system representation of the form
\begin{equation*}
M(q)\ddot{q}+C(q,\dot{q})=Gu\,,
\end{equation*}
with $q=[q^1,q^2,q^3]^T$ and $u=[F_1,F_2]^T$. The symmetric matrix $M(q)$ reads as
\begin{equation*}
M(q)=\left[\begin{array}{ccc}
m_{11}&m_{12}+m_h\Phi(q^3)&m_hq^2\partial_{q^3}\Phi(q^3)\\
* &m_{22}+m_h\Phi(q^3)^2&m_hq^2\Phi(q^3)\partial_{q^3}\Phi(q^3)\\
* & * & m_h+m_h(q^2\partial_{q^3}\Phi(q^3))^2
\end{array}\right]
\end{equation*}
while $C(q,\dot{q})$ as well as $G$ are given by
\begin{equation*}
C(q,\dot{q})=\left[\begin{array}{c}C_1\\q^2EI\smallint_{0}^{L}(\partial^2_z\Phi(z))^2\mathrm{d}z+\Phi(q^3)C_1\\C_1q^2\partial_{q^3}\Phi(q^3)+m_hg
\end{array}\right]\,,\phantom{a}G=\left[\begin{array}{cc}1&0\\0&0\\0&1
\end{array}\right]
\end{equation*}
with \begin{equation*}
C_1=m_h(\dot{q}^3)^2q^2\partial^2_{q^3}\Phi(q^3)+2m_h\dot{q}^3q^2\partial_{q^3}\Phi(q^3)\,.
\end{equation*}
Note that we adopted the abbreviations $m_{11},m_{12}$ and $m_{22}$ from \cite{RamsSchoeberlSchlacher2018}.
If we introduce the time-derivatives of $q$ as states $v$ (velocities), then we get a state representation of the form
\begin{equation}
\label{eq:state_representation_stacker_crane}
\underbrace{\left[\begin{array}{c}
	\dot{q}\\
	\dot{v}
	\end{array}\right]}_{\dot{x}}=\underbrace{\left[\begin{array}{c}
	v\\-M^{-1}(q)C(q,v)+M^{-1}(q)Gu
	\end{array}\right]}_{f(x,u)}\,,
\end{equation}
with $x=[x^1,x^2,x^3,x^4,x^5,x^6]^T=[q^1,q^2,q^3,v^1,v^2,v^3]^T$ and $u=[u^1,u^2]^T=[F_1,F_2]^T$. 
\subsection{Flatness of the Continuous-time Model}
\label{subsec:flatness_analysis_cont}
First, we want to emphasize that the system \eqref{eq:state_representation_stacker_crane} is not static feedback-linearizable. This can be shown with the geometric necessary and sufficient conditions stated e.g. in \cite{NijmeijervanderSchaft:1990}. Nevertheless, as mentioned in \cite{RamsSchoeberlSchlacher2018}, the system \eqref{eq:state_representation_stacker_crane} is flat and possesses a configuration-flat output of the form
\begin{equation}
\label{eq:configuration_flat_output}
y=(x^1+\Phi(x^3)x^2,x^1(m_{22}-m_{12}\Phi(x^3)))\,.
\end{equation}
The flat output \eqref{eq:configuration_flat_output} can be derived systematically using e.g. the approach of \cite{RathinamMurray1996} for underactuated mechanical systems or the more recent geometric methods of  \cite{GstoettnerKolarSchoeberl2021-2}. However, for practical purposes, the configuration-flat output \eqref{eq:configuration_flat_output} is not suitable, because the corresponding parameterizing map \eqref{eq:parameterizing_map_continuous} has singularities at stationary points.
Fortunately, there exists another flat output that can be constructed intuitively using the controller canonical form for linear time-variant systems as shown in \cite{StaudeckerSchlacherHansl2008}. First, by means of the input transformation 
\begin{equation}
\label{eq:input_transformation_stacker_crane}
\bar{u}^1=f^4(x,u)\,,\phantom{aa}\bar{u}^2=f^6(x,u)\,,
\end{equation}
the system \eqref{eq:state_representation_stacker_crane} can be decoupled into two subsystems
\begin{equation}
\label{eq:linear_representation_stacker_crane}
\left[\begin{array}{c}
\dot{x}_1\\\dot{x}_2
\end{array}\right]=\left[\begin{array}{cc}
A_1(x^3,x^6,\bar{u}^2)&0\\
0&A_2
\end{array}\right]\left[\begin{array}{c}
x_1\\x_2
\end{array}\right]+\left[\begin{array}{cc}
b_1(x^3)&0\\0&b_2
\end{array}\right]\left[\begin{array}{c}
\bar{u}^1\\\bar{u}^2
\end{array}\right]\,,
\end{equation}
with $x_1=[x^1,x^2,x^4,x^5]^T$, $x_2=[x^3,x^6]^T$ and
\begin{equation*}
A_2=\left[\begin{array}{cc}
0&1\\0&0
\end{array}\right]\,,\phantom{a}b_2=\left[\begin{array}{c}
0\\1
\end{array}\right]\,.
\end{equation*}
Since the subsystem related to $x_2$ is in controller canonical form, the state $x^3$ (vertical position of the lifting unit) is a candidate for the first component of the flat output $y$. With the choice $y^1=x^3$ we can already parameterize the states $x_2$ and the input $\bar{u}^2$ of the second subsystem by the flat output and its time-derivatives, i.e.
\begin{equation}
\label{eq:parameterizing_map_x3x6u2q}
\begin{aligned}
x^3=y^1\,,\phantom{aa}x^6=\dot{y}^1\,,\phantom{aa}\bar{u}^2=\ddot{y}^1\,.
\end{aligned}
\end{equation}
In order to compute the second component of the flat output, let us fix in an intermediate step a trajectory $y^1(t)$ for the first component of the flat output. The corresponding trajectories for $x^3(t),x^6(t)$ and $\bar{u}^2(t)$ are then fixed via the relations \eqref{eq:parameterizing_map_x3x6u2q}. If we insert $x^3(t),x^6(t)$ and $\bar{u}^2(t)$ into the first subsystem, we get a linear time-variant system
\begin{equation*}
\label{eq:subsystem_time_variant_continuous}
\dot{x}_1=A_1(x^3(t),x^6(t),\bar{u}^2(t))x_1+b_1(x^3(t))\bar{u}^1\,.
\end{equation*}
For linear time-variant flat systems, however, a flat output can be constructed systematically using the state transformation
\begin{equation}
\label{eq:linear_time_variant_state_transformationxx}
\bar{x}_1=T(x^3(t),x^6(t),\bar{u}^2(t),\dot{\bar{u}}^2(t),\dots,{\bar{u}}^{2,(4)}(t))x_1
\end{equation} into the time-variant controller-canonical form, see e.g. \cite{IlchmannMueller2007}.  In such coordinates, the first component of $\bar{x}_1$ corresponds to the second component of the flat output. With the inverse transformation of \eqref{eq:input_transformation_stacker_crane} and \eqref{eq:linear_time_variant_state_transformationxx}, the second component of the flat output can be expressed in terms of the original coordinates and is of the form $y^2=\varphi^2(x,u,\dot{u})$. In contrast to \eqref{eq:configuration_flat_output}, the corresponding parameterizing map for the flat output 
\begin{equation*}
y=(x^3,\varphi^2(x,{u},\dot{{u}}))
\end{equation*} has no singularities at stationary points and can therefore be used for practical purposes. 
\subsection{Motivation and Problem Statement}
The concept of continuous-time flatness allows the systematic design of tracking controllers. Nevertheless, for the practical implementation it is necessary to evaluate the continuous-time control law at a sufficiently high sampling rate. For lower sampling rates, a discrete evaluation of such a continuous-time control law could lead to unsatisfactory results, as we have shown e.g. in \cite{DiwoldKolarSchoeberl2022-1} by means of the laboratory setup of a gantry crane. An appropriate alternative is to design the controller directly for a suitable discretization. The main challenge in the design of discrete-time flatness-based control laws, however, is to find a flat discretization, since in general flatness is not preserved by an approximate or exact discretization of a continuous-time system. For this reason, in \cite{DiwoldKolarSchoeberl2022-1} we proposed a method that combines the explicit Euler-discretization and a prior state transformation in such a way that not only the flatness but also the flat output is preserved. Although this method is restricted to a subclass of continuous-time $x$-flat systems, it can be applied to many practical examples such as the induction motor, the VTOL aircraft, or the unmanned aerial vehicle discussed in \cite{Greeff2021}. Even though the model of a stacker crane \eqref{eq:state_representation_stacker_crane} together with the flat output \eqref{eq:configuration_flat_output} also belongs to this class, in this case our method is not useful for practical applications. The reason is that not only the flat output \eqref{eq:configuration_flat_output} would be preserved but also the corresponding singularities in the parameterizing map. Thus, the main objective of this paper is to construct a flat output for the explicit Euler-discretization of \eqref{eq:state_representation_stacker_crane} in a similar way as it is shown in Section \ref{subsec:flatness_analysis_cont} for the continuous-time model. For this purpose, we subsequently extend the definition of flatness as proposed in \cite{DiwoldKolarSchoeberl2022} to time-variant systems and briefly discuss the controller canonical form for linear time-variant discrete-time systems. Based on these results, in Section \ref{sec:sampled_data_model_stacker_crane_} we assume again that the vertical position of the lifting unit corresponds to the first component of the flat output, and construct the second component by transforming a linear time-variant discrete-time subsystem into the controller canonical form. As we will see below, the constructed flat output depends also on backward-shifts of system variables.
\section{Flatness of Time-variant Discrete-time Systems}
\label{sec:flatness_time_variant}
In the following we recall and extend the definition of flatness as proposed in \cite{DiwoldKolarSchoeberl2022} to time-variant discrete-time systems. 
Note that in this paper we use a slightly different notation than in previous contributions (e.g. \cite{KolarSchoberlDiwold:2019}, \cite{DiwoldKolarSchoeberl2022} and \cite{DiwoldKolarSchoeberl2022-1}), since we focused only on time-invariant discrete-time systems so far.
\subsection{Flatness of Nonlinear Time-variant Systems}
\label{sec:flatness_time_variant_nonlinear_sys}
\label{sec:Flatness-of-Discrete-time}
We consider time-variant discrete-time systems in state representation\footnote{Instead of $f^{i}(k,x_k,u_k)$ we write $f_k^{i}(x_k,u_k)$ for a better readability.}
\begin{equation}
\label{eq:system_equations_disc}
x_{k+1}=f_k(x_k,u_k)
\end{equation}
with $\dim(x_k)=n$, $\dim(u_k)=m$ and smooth functions $f_k^{i}(x_k,u_k)$.
Furthermore, we assume that the Jacobian matrix of
$f_k$ w.r.t. $(x_k,u_k)$ meets
\begin{equation}
\textrm{rank}(\partial_{(x_k,u_k)}f_k)=n
\label{eq:submersivity_condition}
\end{equation}
for all $k$ (submersivity). 
As stated in \cite{DiwoldKolarSchoeberl2022}, flatness for discrete-time systems can be characterized by a one-to-one correspondence between the system trajectories $(x_k,u_k)$ and trajectories $y_k$ of a trivial system, i.e., arbitrary trajectories that need
not satisfy any difference equation. In contrast to the time-invariant systems considered in \cite{DiwoldKolarSchoeberl2022}, for time-variant systems \eqref{eq:system_equations_disc} the maps
\begin{equation}
\label{eq:parm_map_time_variant}
(x_k,u_k)=F(k,y_{k-r_1},\dots,y_k,\dots,y_{k+r_2})
\end{equation}
and
\begin{equation}
\label{eq:flat_output_time_variant}
\begin{aligned}
y_k=\varphi(k,x_{k-q_1},u_{k-q_1},\dots,x_k,u_k,\dots,x_{k+q_2},u_{k+q_2})\
\end{aligned}
\end{equation}
describing this one-to-one correspondence may also depend explicitly on $k$. As explained in \cite{DiwoldKolarSchoeberl2022}, the composition of \eqref{eq:parm_map_time_variant} with the occurring shifts of \eqref{eq:flat_output_time_variant} must yield the identity map, and vice versa. Furthermore, since the trajectories $y_k$ of a trivial system are arbitrary, substituting \eqref{eq:parm_map_time_variant} into \eqref{eq:system_equations_disc} must also yield an identity. The map \eqref{eq:flat_output_time_variant} can be further simplified by taking into account that every trajectory $\dots,x_{k-1},u_{k-1},x_k,u_k,x_{k+1},u_{k+1},\dots$ satisfies \eqref{eq:system_equations_disc}. Hence, the forward-shifts of $x_k$ can be expressed as functions of $k,x_k,u_k,u_{k+1},\dots$ via the successive compositions
\begin{equation}
\label{eq:to_replace_forward_shifts_of_x}
\begin{aligned}
x_{k+1}&=f_k(x_k,u_k)\\
x_{k+2}&=f_{k+1}(f_k(x_k,u_k),u_{k+1})\,.\\
\vdots
\end{aligned}
\end{equation}
A similar argument holds for the backward-direction. Since \eqref{eq:system_equations_disc} meets the submersivity condition \eqref{eq:submersivity_condition}, there always exist $m$ functions $g_k(x_k,u_k)$ such that the map
\begin{equation}
\label{eq:extending_f_to_diffeomorphism}
x_{k+1}  =f_k(x_k,u_k)\,,\phantom{aa}\zeta_k=g_k(x_k,u_k)\,,
\end{equation}
is locally invertible.
If we denote by $(x_k,u_k)=\psi_k(x_{k+1},\zeta_k)$ its inverse
\begin{align}
\label{eq:inverted_extension}
\begin{aligned}x_k & =\psi_{x,k}(x_{k+1},\zeta_k)\,,\end{aligned}
& u_k=\psi_{u,k}(x_{k+1},\zeta_k)\,
\end{align}
then we can express the backward-shifts of $(x_k,u_k)$ as functions of $k,x_k,\zeta_{k-1},\zeta_{k-2},\dots$ via the successive compositions
\begin{equation}
\label{eq:to_replace_backward_shifts_of_x_and_u}
\begin{aligned}
(x_{k-1},u_{k-1})&=\psi_{k-1}(x_k,\zeta_{k-1})\\
(x_{k-2},u_{k-2})&=\psi_{k-2}(\psi_{x,k-1}(x_k,\zeta_{k-1}),\zeta_{k-2})\,.\\
\vdots
\end{aligned}
\end{equation}
Thus, with \eqref{eq:to_replace_forward_shifts_of_x} and \eqref{eq:to_replace_backward_shifts_of_x_and_u}, the map \eqref{eq:flat_output_time_variant} can be written as
\begin{equation}
\label{eq:flat_output_as_functions_of_zeta_and_u}
y_k=\varphi(k,\zeta_{k-q_1},\dots,\zeta_{k-1},x_k,u_k,\dots,u_{k+q_2})\,.
\end{equation}
The map \eqref{eq:flat_output_as_functions_of_zeta_and_u} is called a flat output of the time-variant system \eqref{eq:system_equations_disc}, and the corresponding map \eqref{eq:parm_map_time_variant} describes the parameterization of the system variables $(x_k,u_k)$ by this flat output.
\begin{rem}
	\label{rem:choice_zeta}The flatness
	of the system (\ref{eq:system_equations_disc}) does not depend on the choice
	of the functions $g_k(x_k,u_k)$ of \eqref{eq:extending_f_to_diffeomorphism}. The representation \eqref{eq:flat_output_as_functions_of_zeta_and_u}
	of the flat output may differ, but the parameterization \eqref{eq:parm_map_time_variant}
 is not affected. 
\end{rem}
Note also that considering backward-shifts in both \eqref{eq:parm_map_time_variant} and \eqref{eq:flat_output_as_functions_of_zeta_and_u} is actually not necessary. If there exist a parameterizing map \eqref{eq:parm_map_time_variant} and a flat output \eqref{eq:flat_output_as_functions_of_zeta_and_u}, then one can always define a new flat output as the $r_1$-th backward-shift of the original flat output. The corresponding parameterizing map is then of the form
\begin{equation}
\label{eq:parameterizing_map_forward_shifts_only}
(x_k,u_k)=F(k,y_k,\dots,y_{k+r})
\end{equation}
with $r=r_1+r_2$ and does not contain backward-shifts. Systems that possess a flat output  $y_k=\varphi(k,x_k,u_k,\dots,u_{k+q_2})$ which is independent of backward-shifts of $\zeta_k$ together with a parameterizing map \eqref{eq:parameterizing_map_forward_shifts_only} will be denoted in the following as forward-flat.
\subsection{Flatness of Linear Time-variant Single-input Systems}
In this section, we focus on the special case of linear time-variant discrete-time systems\footnote{Again, instead of $A(k)$ and $b(k)$, we write $A_k$ and $b_k$, respectively.}
\begin{equation}
\label{eq:linar_time_variant_discrete_time_system}
x_{k+1}=A_{k}x_k+b_ku_k
\end{equation} 
with $\dim(x_k)=n$ and $\dim(u_k)=1$. For constructing a flat output of the sampled-data model of the single mast stacker crane, like in the continuous-time case a linear time-variant transformation
\begin{equation}
\label{eq:linear_time_variant_state_transformation}
x_k=T_k\bar{x}_k
\end{equation}
into the controller canonical form is useful. In general, the transformation \eqref{eq:linear_time_variant_state_transformation} yields again a linear time-variant system of the form
\begin{equation}
\label{eq:transformed_linear_time_variant_system_bar_notation}
\bar{x}_{k+1}=\bar{A}_k\bar{x}_{k}+\bar{b}_ku_k
\end{equation} 
with $\bar{A}_k=T^{-1}_{k+1}A_kT_k$ and $\bar{b}_k=T^{-1}_{k+1}b_k$.
A system \eqref{eq:transformed_linear_time_variant_system_bar_notation} is said to be in controller canonical form if the matrix $\bar{A}_k$ and the vector $\bar{b}_k$ read as
\begin{equation*}
\label{eq:matrices_of_controller_canonical_form}
\bar{A}_k=\left[\begin{array}{cccc}
0 & 1 & \cdots & 0\\
\vdots & \vdots & \ddots & \vdots\\
0 & 0 & \cdots & 1\\
-a_{0,k} & -a_{1,k} & \cdots & -a_{n-1,k}
\end{array}\right]\,,\phantom{a}\bar{b}_k=\left[\begin{array}{c}0\\\vdots\\0\\1
\end{array}\right]
\end{equation*}
with time-variant coefficients $a_{0,k},\dots,a_{n-1,k}$, see also \cite{HwangKamenTeply1985}. Note that $a_{0,k},\dots,a_{n-1,k}$ do not necessarily correspond to the coefficients of the characteristic polynomial of $A_k$. Obviously, the existence of a transformation into controller canonical form is a sufficient condition for flatness, since the first component of $\bar{x}_k$ forms a flat output. In the following, we show how the transformation into controller canonical form can be constructed in a systematic way. For this purpose, we use an ansatz $y_k=c^T_k x_k$
and require that all forward-shifts of $y_k$ up to the order $n-1$ are independent of $u_k$, i.e.,
\begin{equation}
\label{eq:forward_shifts_of_ansatz_up_to_n-1}
\begin{aligned}
y_{k+1}&=c^T_{k+1}A_kx_k+\overbrace{c^T_{k+1}b_k}^{=0}u_k\\
y_{k+2}&=c^T_{k+2}A_{k+1}A_kx_k+\overbrace{c^T_{k+2}A_{k+1}b_k}^{=0}u_k\,.\\
&\phantom{a}\vdots
\end{aligned}
\end{equation}
For the $n$-th forward-shift of $y_k$ we require
\begin{equation}
\label{eq:forward_shifts_of_ansatz_n}
\begin{aligned}
y_{k+n}=&c^T_{k+n}A_{k+n-1}\cdots A_kx_k+\overbrace{c^T_{k+n}A_{k+n-1}\cdots A_{k+1}b_k}^{=1}u_k\,.
\end{aligned}
\end{equation}
By applying suitable backward-shifts, the conditions \eqref{eq:forward_shifts_of_ansatz_up_to_n-1} and \eqref{eq:forward_shifts_of_ansatz_n} can be written as
\begin{equation}
\label{eq:equation_for_computation_of_c_k}
c^T_k\underbrace{[
	b_{k-1},A_{k-1}b_{k-2},\dots,A_{k-1}\cdots A_{k-n+1}b_{k-n}
	]}_{:=M_k}=e_n^T
\end{equation}
with the $n$-th unit vector $e_n$. Obviously, there exists a solution for $c_k^T$ if and only if the matrix\footnote{For linear time-invariant systems the matrix $M_k=M$ corresponds to the controllability matrix $M=[b,Ab,\dots,A^{n-1}b]$.} $M_k$ is regular for all $k$. The transformation into controller canonical form is then given by
\begin{equation}
\label{eq:transformation_into_controller_canonical_form}
\bar{x}_k=\left[\begin{array}{c}
y_k\\y_{k+1}\\\vdots\\y_{k+n-1}
\end{array}\right]
=\underbrace{\left[\begin{array}{c}
	c^T_k\\c^T_{k+1}A_k\\\vdots\\c^T_{k+n-1}A_{k+n-2}\cdots A_k
	\end{array}\right]}_{T_k^{-1}}x_k\,.
\end{equation}
In the controller canonical form \eqref{eq:transformed_linear_time_variant_system_bar_notation}, the first component of the state $\bar{x}_k$ is a flat output. For the system \eqref{eq:linar_time_variant_discrete_time_system} in original coordinates, this flat output is given by
\begin{equation}
\label{eq:flat_output_linear_timevariant}
y_k=\varphi(k,x_k)=c^T_kx_k\,.
\end{equation}
The corresponding parameterization of $(x_k,u_k)$ follows as
\begin{equation}
\label{eq:parameterization_of_x_via_controller_canonical_form}
x_k=F_x(k,y_k,\dots,y_{k+n-1})=T_k\bar{x}_k
\end{equation}
and
\begin{equation}
\label{eq:parameterization_of_u_via_controller_canonical_form}
u_k=F_u(k,y_k,\dots,y_{k+n})=y_{k+n}+\sum_{i=0}^{n-1}a_{i,k}y_{k+i}\,.
\end{equation}
Since neither \eqref{eq:flat_output_linear_timevariant} nor \eqref{eq:parameterization_of_x_via_controller_canonical_form} and \eqref{eq:parameterization_of_u_via_controller_canonical_form} depend on backward-shifts, the condition
 \begin{equation}
	\label{eq:reachability_rank_condition}
	\mathrm{rank}(M_k)=n\,,\,\forall k
	\end{equation} 
	even ensures forward-flatness. In fact, it can be shown that this condition is not only sufficient but also necessary for flatness. 
\section{Sampled-data Model of the Single Mast Stacker Crane}
\label{sec:sampled_data_model_stacker_crane_}
In this section we consider the explicit Euler-discretization
\begin{equation}
\label{eq:explicit_euler_disc}
x_{k+1}=x_k+T_sf(x_k,u_k)
\end{equation}
of the single mast stacker crane \eqref{eq:state_representation_stacker_crane} and derive a flat output similar as in the continuous-time case. Subsequently, we use the corresponding parameterizing map to compute a feedforward control which achieves a transition of the single mast stacker crane between two rest positions.
\subsection{Flatness of the Sampled-data Model}
\label{sec:sampled_data_model_stacker_crane}
In the following, we proceed with the explicit Euler-discretization \eqref{eq:explicit_euler_disc} of the transformed system \eqref{eq:linear_representation_stacker_crane}, i.e.
\begin{equation*}
x_{k+1}=x_k+T_sf(x_k,\bar{u}_k)\,.
\end{equation*}
This is possible since the explicit Euler-discretization and input-transformations commute\footnote{If two continuous-time systems are related by an input transformation $\bar{u}=\Phi_u(x,u)$, then their explicit Euler-discretizations are related by the same input transformation $\bar{u}_k=\Phi_u(x_k,u_k)$.}, see e.g. \cite{DiwoldKolarSchoeberl2022-1}. The Euler-discretized system reads as
\begin{equation*}
\left[\begin{array}{c}
x_{1,k+1}\\x_{2,k+1}
\end{array}\right]=\left[\begin{array}{cc}
A_1(x^3_k,x^6_k,\bar{u}^2_k)&0\\
0&A_2
\end{array}\right]\left[\begin{array}{c}
x_{1,k}\\x_{2,k}
\end{array}\right]+\left[\begin{array}{cc}
b_1(x^3_k)&0\\0&b_2
\end{array}\right]\bar{u}_k
\end{equation*}
with $x_{1,k}=[x^1_k,x^2_k,x^4_k,x^5_k]^T$, $x_{2,k}=[x^3_k,x^6_k]^T$ and
\begin{equation}
\label{eq:structurally_flat_triangular_form_subsystem_2_disc}
A_2=\left[\begin{array}{cc}
1&T_s\\0&1
\end{array}\right]\,,\phantom{a}b_2=\left[\begin{array}{c}
0\\T_s
\end{array}\right]\,.
\end{equation}
Note that also the matrix $A_1$ and the vector $b_1$ differ from the ones defined in \eqref{eq:linear_representation_stacker_crane}. Like in the continuous-time case, the vertical position of the lifting unit $x^3_k$ is a candidate for $y^1_k$, since with \eqref{eq:structurally_flat_triangular_form_subsystem_2_disc} we can already parameterize
\begin{equation}
\label{eq:parameterization_x3_x6_u2q}
\begin{aligned}
x^3_k&=y^1_k\,,\phantom{a}x^6_k=\tfrac{y^1_{k+1}-y^1_{k}}{T_s}\,,\phantom{a}\bar{u}^2_k=\tfrac{y^1_{k+2}-2y^1_{k+1}+y^1_k}{T_s^2}\,.
\end{aligned}
\end{equation}
Once a trajectory $y^1_k$ is fixed, there remains a linear time-variant system of the form
\begin{equation*}
\label{eq:linear_time_variant_disc_subsystem_1}
x_{1,k+1}=A_1(x^3_k,x^6_{k},\bar{u}^2_k)x_{1,k}+b_1(x^3_k)\bar{u}^1_k\,,
\end{equation*}
which can be checked w.r.t. flatness by condition \eqref{eq:reachability_rank_condition}. Since there appear backward-shifts of $A_1(x^3_k,x^6_{k},\bar{u}^2_k)$ and $b_1(x^3_k)$ within the matrix $M_k$, we need to introduce the functions $g_k(x_k,u_k)$ of \eqref{eq:extending_f_to_diffeomorphism}. With the inverse \eqref{eq:inverted_extension}, we can then express the backward-shifts of $x^3_k$ as functions of $x_k$ and backward-shifts of $\zeta_k$, see also equation \eqref{eq:to_replace_backward_shifts_of_x_and_u}. With the choice
\begin{equation*}
\zeta^1_k=x^3_k\,,\phantom{aaa}\zeta^2_k=x^1_k\,,
\end{equation*}
we get the matrix
\begin{equation*}
\label{eq:reachability_matrix_stacker_crane}
\begin{aligned}
M_k=[&b_1(\zeta^1_{k-1}),A_1(\zeta^1_{k-1},x^3_{k},x^6_{k})b_1(\zeta^1_{k-2}),\dots,\\
&A_1(\zeta^1_{k-1},x^3_{k},x^6_{k})\cdots A_1(\zeta^1_{k-3},\zeta^1_{k-2},\zeta^1_{k-1})b_1(\zeta^1_{k-4})]
\end{aligned}
\end{equation*}
 which is regular for arbitrary values $\zeta^1_{k-4},\dots,\zeta^1_{k-1},x^3_k,x^6_k$. This holds independently of the chosen extension $\zeta_k=g_k(x_k,u_k)$. With the corresponding controller canonical form, we can now determine the second component of the flat output. Due to relation \eqref{eq:equation_for_computation_of_c_k}, also $c^T_k$ depends on $\zeta^1_{k-4},\dots,\zeta^1_{k-1},x^3_k,x^6_k$, and the second component of the flat output $y_k$ has the form
 \begin{equation*}
 y^2_k=c_k^Tx_1=c^T(\zeta^1_{k-4},\dots,\zeta^1_{k-1},x^3_k,x^6_k)x_1\,.
 \end{equation*}
In order to obtain the parameterization for $x_1$ and $\bar{u}^1$, we use the matrix $T_k$ of \eqref{eq:transformation_into_controller_canonical_form}, which is of the form
\begin{equation*}
\label{eq:T_as_function_of_zeta_x_u}
T_k=T(\zeta^1_{k-4},\dots,\zeta^1_{k-1},x^3_k,x^6_k,\bar{u}^2_k,\dots,\bar{u}^2_{k+2})\,.
\end{equation*}
With the parameterization \eqref{eq:parameterization_x3_x6_u2q} as well as its forward- and backward-shifts, the matrix $T_k$ can be expressed in terms of $y^1_k$ and its shifts as $T_k=T(y^1_{k-4},\dots,y^1_{k+4})$. Evaluating \eqref{eq:parameterization_of_x_via_controller_canonical_form} and \eqref{eq:parameterization_of_u_via_controller_canonical_form} yields a parameterization of the form
\begin{equation*}
(x_k,\bar{u}_k)=F(y^1_{k-4},\dots,y^1_{k+5},y^2_k,\dots,y^2_{k+4})
\end{equation*}
together with the flat output
\begin{equation*}
y_k=(x^3_k,\varphi^2(\zeta^1_{k-4},\dots,\zeta^1_{k-1},x_k))\,.
\end{equation*}
As mentioned above, considering backward-shifts both in the parameterizing map as well as the flat output can be avoided by redefining $y^1_{k-4}\rightarrow y^1_{k}$, which leads to
\begin{equation}
\label{eq:parameterizing_map_stacker_crane}
(x_k,\bar{u}_k)=F(y^1_{k},\dots,y^1_{k+9},y^2_k,\dots,y^2_{k+4})\,
\end{equation}
together with
\begin{equation*}
\label{eq:flat_output_disc_stacker_crane}
y_k=(\zeta^1_{k-4},\varphi^2(\zeta^1_{k-4},\dots,\zeta^1_{k-1},x_k))\,.
\end{equation*}
Since input transformations and the Euler-discretization commute, we can also parameterize the original input $u_k$ of \eqref{eq:explicit_euler_disc} by means of the inverse of \eqref{eq:input_transformation_stacker_crane}. Thus, we can finally conclude that the Euler-discretization \eqref{eq:explicit_euler_disc} of the original system \eqref{eq:state_representation_stacker_crane} is also flat.
\begin{rem}
Since the computed flat output of the sampled-data system \eqref{eq:explicit_euler_disc} depends on backward-shifts of the system variables, the question arises whether there also exists a forward-flat output. By using the geometric necessary and sufficient conditions for forward-flatness derived in \cite{KolarSchoberlDiwold:2019} and \cite{KolarDiwoldSchoberl:2019}, it can be shown that this is not the case. Thus, the sampled-data model of the stacker crane is indeed a practical example for which the usual definition of flatness (without taking backward-shifts into account) is not sufficient.
\end{rem}

\subsection{Trajectory Planning and Feedforward Control}
\label{subsec:trajectory_planning}
\begin{figure}
	\begin{center}
		\includegraphics[width=8.4cm]{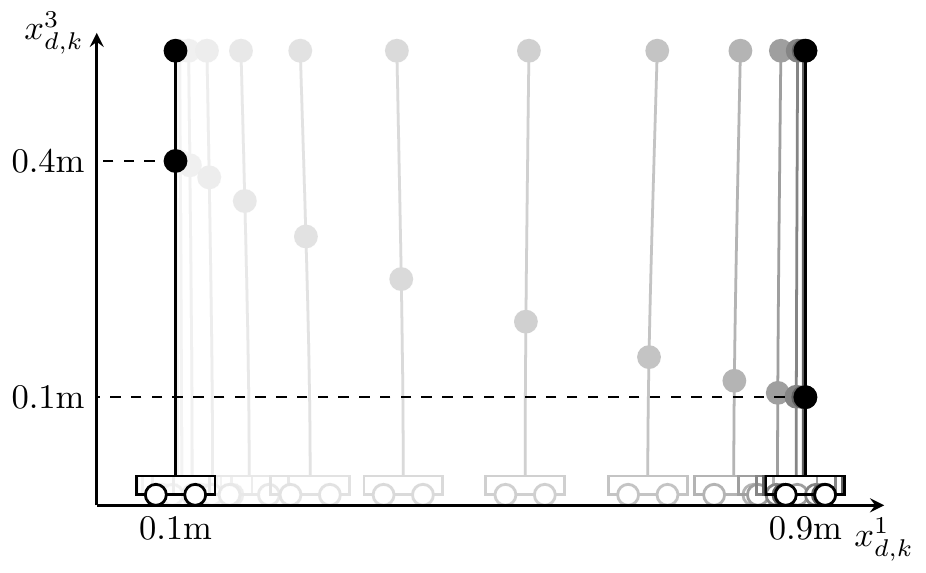}    
		\caption{Transition between two rest positions.} 
		\label{fig:movement}
	\end{center}
\end{figure}

\begin{figure}
	\begin{center}
		\includegraphics[width=7.7cm]{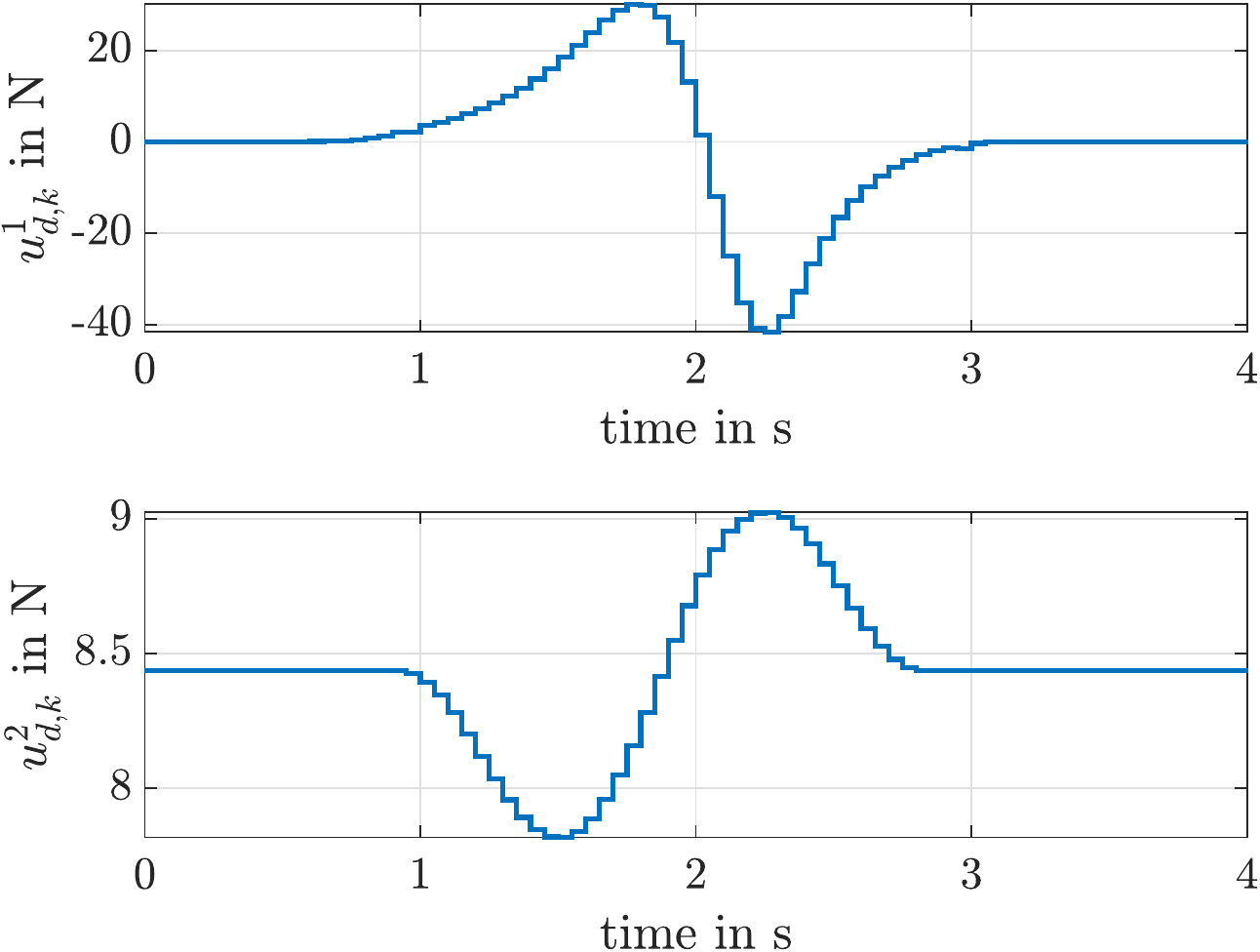}    
		\caption{Feedforward control ($T_s=50\text{ms}$).} 
		\label{fig:feedforward}
	\end{center}
\end{figure}

The derived parameterizing map can be used efficiently for planning trajectories. We illustrate this by planning a trajectory which achieves a transition between two rest positions $x_0$ at $k=0$ and $x_N$ at $k=N$. In rest positions, the shape of the beam is straight and the velocities are zero. Hence, the corresponding state variables meet $x^i_0=x^i_N=0$ for $i=2,4,5,6$. Furthermore, since in a rest position the shifts of the flat output are constant, a desired reference trajectory must meet
\begin{equation}
\label{eq:setting_y_d_constant}
\begin{aligned}
y^1_{d,0}&=\cdots=y^1_{d,9}\,,&\phantom{a}y^2_{d,0}&=\cdots=y^2_{d,4}\\
y^1_{d,N}&=\cdots=y^1_{d,N+9}\,,&\phantom{a}y^2_{d,N}&=\cdots=y^2_{d,N+4}\,.\\
\end{aligned}
\end{equation} 
If we insert \eqref{eq:setting_y_d_constant} into the parameterization \eqref{eq:parameterizing_map_stacker_crane} for $x_k$, we get
\begin{equation*}
\begin{aligned}
x_0&=F_x(y^1_{d,0},\dots,y^1_{d,0},y^2_{d,0},\dots,y^2_{d,0})\\
x_N&=F_x(y^1_{d,N},\dots,y^1_{d,N},y^2_{d,N},\dots,y^2_{d,N})\,,
\end{aligned}
\end{equation*}
which can be solved for $y^1_{d,0},y^2_{d,0}$ and $y^1_{d,N},y^2_{d,N}$. The remaining values of $y^1_{d,10},\dots,y^1_{d,N-1}$ and $y^2_{d,5},\dots,y^2_{d,N-1}$ can be chosen arbitrarily. The simplest approach is to use polynomial functions $y^i_{d,k}=\sum_{l=0}^{l_i}p^i_l k^l$ of appropriate order $l_i$ with suitable coefficients $p^i_l$. Once the trajectories $y^1_{d,k}$ and $y^2_{d,k}$ are fixed, the corresponding state- and input-trajectories  $(x_{d,k},u_{d,k})$ are uniquely determined by the parameterizing map. In Fig. \ref{fig:movement} the transition between two rest positions is illustrated, and the corresponding feedforward control is depicted in Fig. \ref{fig:feedforward}. All parameters for the model \eqref{eq:state_representation_stacker_crane} as well as the Rayleigh-Ritz ansatz function $\Phi(z)$ were adopted from \cite{StaudeckerSchlacherHansl2008} and \cite{RamsSchoeberlSchlacher2018}. The sampling time was set to $T_s=50\text{ms}$.

\section{Conclusion and Outlook}
In this contribution, we have shown that the explicit Euler-discretization of the single mast stacker crane \eqref{eq:state_representation_stacker_crane} is flat. For this purpose, we have constructed a flat output by means of the controller canonical form for linear time-variant discrete-time systems. In contrast to the sampled-data model of a gantry crane that we derived in \cite{DiwoldKolarSchoeberl2022-1}, the explicit Euler-discretization of \eqref{eq:state_representation_stacker_crane} is flat but not forward-flat. However, as illustrated in Section \ref{subsec:trajectory_planning}, also the more general class of flat systems (including backward-shifts) allows for the systematic planning of trajectories and design of feedforward controls. Like in the continuous-time case, one could also formulate an optimization problem (e.g. time-optimal transition) instead of using polynomial functions for $y_{d,k}$. While often ansatz functions (e.g. spline curves) are used for this task in the continuous-time case in order to obtain a finite number of optimization variables, the optimization problem is inherently finite-dimensional in the discrete-time case. Apart from formulating and solving an optimization problem, subject to further research could be the design of discrete-time flatness-based tracking controllers based on an exact linearization as proposed in \cite{KolarDiwoldGstoettnerSchoeberl2022}. Subsequently, the discrete-time approach could be compared to existing continuous-time (sampled) control laws in simulations as well as on the laboratory setup of the single mast stacker crane.



                                                   
\bibliography{Bibliography_Johannes_Mai_2022}             

\end{document}